\documentstyle[12pt]{article}
\newcommand{\be}{\begin{equation}}
      \newcommand{\ee}{\end{equation}}
      \newcommand{\ba}{\begin{eqnarray}}
       \newcommand{\ea}{\end{eqnarray}}
\newcommand{\ban}{\begin{eqnarray*}}
       \newcommand{\ean}{\end{eqnarray*}}

\newcommand{\ra}{\rightarrow}

 \newcommand{\qed}{\hspace*{\fill}\rule{3mm}{3mm}\quad}
 \newcommand{\Pf}{\noindent {\em Proof.} }

\newcommand{\sect}[1]{\section{#1} \setcounter{equation}{0}}

\newtheorem{theo}{Theorem}[section]

\begin{document}
\newtheorem{lem}[theo]{Lemma}
\newtheorem{prop}[theo]{Proposition}  
\newtheorem{coro}[theo]{Corollary}

\title{Smoothing Riemannian Metrics with Ricci Curvature Bounds}
\author{Xianzhe Dai\thanks{Partially supported by NSF Grant \# DMS9204267 and Alfred  P. Sloan Fellowship} \\
  \and Guofang Wei\thanks {Partially supported by NSF Grant \# DMS9409166.
}\\
\and Rugang Ye\thanks{Partially supported by NSF Grant \# DMS9401106}
}
\date{}
\maketitle

\begin{abstract}
We prove that Riemannian metrics with an absolute Ricci curvature bound and a conjugate
radius bound can be smoothed to having a sectional curvature bound. Using this
we derive a number of results about structures of manifolds with Ricci
curvature bounds.
\end{abstract}

\newcommand{\Tr}{\mbox{Tr}}
\newcommand{\inj}{\mbox{inj}}
\newcommand{\vol}{\mbox{vol}}
\newcommand{\diam}{\mbox{diam}}
\newcommand{\Ric}{\mbox{Ric}}
\newcommand{\conj}{\mbox{conj}}
\newcommand{\Hess}{\mbox{Hess}} 
\newcommand{\divg}{\mbox{div}}
\sect{Introduction}

A central theme in global Riemannian geometry is to understand the global geometric and topological structures of manifolds with appropriate curvature bounds. In this regard, manifolds with sectional curvature bounds have been
better understood than those with Ricci curvature bounds. It is natural to try
to apply the results and methods in the former case to the latter. Hence we ask 
the following question: can one deform or ``smooth" a metric with a Ricci
curvature bound to a metric with a sectional curvature bound?
In this paper we would like to address this question.

Previous work on smoothing mainly deals with metrics having already some sort of
sectional curvature bounds. For example, one can smooth a metric with a 
sectional curvature bound to a metric with bounds on all the covariant
derivatives of its Riemannian curvature tensor. Geometric applications of
this smoothing  result can be found in [CG1] and [CFG]. Its generalizations
to weaker  bounds on sectional curvatures can be found in e.g. [Yn1, Yn2]. So far, two
main methods of smoothing have been applied. The works [BMR], [B], [S] and [Yn1,2] ,
among others, use the Ricci flow, while [CG2] and [Ab] use an embedding method.
The Ricci flow, as a heat flow type equation, tends to average out geometric
quantities, and hence is a natural tool for smoothing. The embedding method
employs an elmentary averaging process in a linear space in which a given
Riemannian manifold is embedded isometrically. 

The following is a precise formulation of the problem of smoothing Riemannian
metrics with Ricci curvature bounds.

\noindent {\bf Question}: Given a positive integer $n$ and any $\epsilon > 0$, does there exist a constant $C(n, \epsilon)$ such that for any closed Riemannian manifold 
$(M^n,\  g)$ with $|\Ric(g)| \leq 1$ we can find another Riemannian metric $\tilde{g}$ on $M$ satisfying \\
1) $ |\tilde{g}-g| <\epsilon$ and \\
2) $ |K(\tilde{g})| \leq C(n, \epsilon)$?

By examples in [A1] among others, the answer to the question in this generality is actually negative. Consequently, additional geometric conditions are needed 
for smoothing metrics with Ricci curvature bounds. In this note we consider
bounds on conjugate radius. \\

\noindent {\bf Definition}. We define for $n \geq 2$ and $r_0 > 0$
\[
{\cal M}(n,,r_0) = \left\{ (M,g) \left| \begin{array}{l} (M,g) \mbox{ is a compact Riemannian manifold}, \\
\dim M = n, \ 
|\Ric | \leq 1\ \mbox{and}\ \conj \geq r_0 
\end{array} \right.
\right\},
\]
where ``conj" denotes the conjugate radius of $M$. \\

Recall that the conjugate radius of a manifold $M$ can be defined to be 
the maximal radius $r$ such that for every $q \in M$, the exponential map
exp$_q$ has maximal rank in the open ball of radius $r$ centered at the origin
of the tangent space $T_qM$. In order to understand the geometric content of
bounds on conjugate radius, it is good to compare with bounds on injectivity
radius. As is well-known, the corresponding space $\tilde{\cal M}(n,r_0)$
concerning injectivity radius (replacing conj by injectivity radius in the above
definition; or equivalently , adding a lower bound on volume \cite{cgt}) is
$C^{1,\alpha}$ precompact \cite{a2}. Hence the situation is quite simple. On the
other hand, under a conjugate radius bound, manifolds can collapse, and
the space ${\cal M}(n,,r_0)$ contains an abundance of geometric and topological
structures. Note also that if $K_M \leq K$ then $\conj \geq
\frac{\pi}{\sqrt{K}}$. This is the usual way of controlling conjugate radius.
But we emphasize that conjugate radius is an important independent geometric invariant.

It seems that Ricci flow is more convenient than the embedding method in our situation, so we choose it. We have
\begin{theo} \label{main}
There exist $T(n,r_0) > 0$ and $C(n,r_0)>0$ such that for any manifold $(M,g_0) \in {\cal M}(n,,r_0)$, the Ricci flow 
\be  \label{rf}
\frac{\partial g}{\partial t} = -2 \Ric (g), \ \ g(0) = g_0
\ee
has a  unique smooth solution $g(t)$ for $0 \leq t \leq T(n,r_0)$ satisfying 
\ba
|g(t) - g_0| & \leq &  4t \label{t1},\\
|Rm (g(t))|_\infty &\leq & C(n,r_0) t^{-1/2}  \label{t2}, \\
|\Ric (g(t))|_\infty & \leq & 2.
\ea
\end{theo}
\noindent {\bf Remark}. Here the conjugate radius condition is only used in the following way: it implies that for every $x \in M$ the lifted metric on $B_x(r_0) \subset T_xM$ has uniform $L^{2,p}$ harmonic coordinates, see Sections 3 and 4.

Using Theorem~\ref{main}, we are able to extend a number of results concerning manifolds with sectional curvature bounds to manifolds with Ricci curvature bounds. First, we have the following generalization of Gromov's almost flat manifold theorem \cite{g1}.

\begin{theo} \label{aflat}
There exists an $\epsilon =\epsilon(n,r_0) > 0$ such that if a manifold $M^n$ has a metric with $|\Ric_M| \leq 1$, conj $\geq r_0$ and diam $\leq \epsilon$, then $M$ is diffeomorphic to an infranilmanifold.
\end{theo}

 Next we have a generalization of Fukaya's fibration theorem \cite{f}.

\begin{theo} \label{fib}
There exists an $\epsilon =\epsilon (n,r_0, \mu) > 0$ such that if $M^n$ and $N^m$ ($m\leq n$) are compact manifolds with $M^n\in {\cal M}(n,r_0)$, $|\Ric_N | \leq 1, \inj_N \geq \mu$, and if the Hausdorff distance $d_H(M,N) \leq \epsilon(n,r_0,\mu)$, then there exists a fibration $f: M \ra N$ such that \\
1) The fibre of $f$ is an almost flat manifold. \\
2) $f$ is an almost Riemannian submersion.
  \end{theo}

Another result is an extension of the double soul theorem in \cite{pz} (it can be used to derive a sphere theorem, see \cite{pz} for details). For a compact manifold $M$ let $ex(M) = \min_{p,q} \max_x (d(p,x)+d(x,q)- d(p,q))$. We have
\begin{theo}  \label{soul}
There exists a positive constant $\epsilon(n,r_0,D)$ such that if $M^n\in {\cal M}(n,r_0)$ with diam$_M \leq D$ and $ex(M) \leq \epsilon$, then either $M$ is infranil or diffeomorphic to the union of two normal bundles over two embedded infranilmanfolds in $M$.
\end{theo}

Let Min\,Vol$(M)$ be the infimum of the volumes of all the complete Riemannian metrics on $M$ with sectional curvature $|K| \leq 1$. We can extend the gap result for minimal volume in dimension four in \cite{ro1}.
\begin{theo}  \label{minv}
There exists a real number $v(r_0) > 0$ such that if a $4$-manifold $M$ has a metric with $|\Ric_M| \leq 1$, conj $\geq r_0$ and Vol $\leq v(r_0)$, then 
 Min\,Vol$(M) = 0$.
\end{theo}

As an immediate corollary of Theorem~\ref{main} and Gromov's uniform betti number estimate \cite{g2}, we have that the betti numbers of manifolds in ${\cal M}(n,r_0)$ are uniformly bounded. Indeed, such a uniform bound still holds if the absolute Ricci curvature bound is replaced by a lower bound \cite{w}. Note however that ${\cal M}(n,r_0)$ contains infinitely many homotopy types.
With Theorem~\ref{fib} at disposal, we can also say something about the higher
homotopy groups, generalizing Rong's results \cite{ro2}. \begin{theo}  \label{h2}
For each $D > 0$ and each $q \geq 2$, the $q$-th rational homotopy group $\pi_q (M)\otimes Q$ has at most $N(n,D,q)$ isomorphism classes for manifolds in ${\cal M}(n,r_0)$ with diam$_M \leq D$.
\end{theo}

Concerning the fundamental group we have
\begin{theo}  \label{fg}
For manifolds in ${\cal M}(n,r_0)$ with diam$_M \leq D$, $\pi_1 (M)$ has a normal nilpotent subgroup $G$ such that \\
1) The minimal number of generators for $G$ is less than or equal to $n$, \\
2) $\pi_1(M)/G$ has at most $N(n,r_0,D)$ isomorphic classes up to a possible normal $Z_2$-extension.
\end{theo} 

Now we would like to say a few words about the strategy of proving Theorem~\ref{main}. The Ricci flow exists at least for a very short time, and the main point is to show that it actually exists on a time interval of uniform size on which the desired estimates hold. A simple, but crucial idea is to lift the solution metrics to suitable domains of the tangent space via the initial exponential maps. The size of these domains can be made uniform by the virtue of the conjugate radius bound. Applying the results in~\cite{a2} we obtain an initial $L^p$ bounds on the sectional curvatures and an initial bound on the Sobolev constant for the lifted metrics. The lifted metrics still satisfy the Ricci flow equation, but they are only defined locally and no a priori control near the boundary is given for them. We apply Moser's weak maximum principle to estimate the sectional curvatures of the lifted metrics in a way similar to [Yn1]. As in [Yn1], we need to control the $L^p$ bound of the sectional curvatures of the lifted metrics along the flow. The subtlety here is how to handle the difficulty caused by lack of control near boundary. Fortunately we found a covering argument to resolve this problem.

It may be possible to obtain part of the smoothing result of Theorem~\ref{main} by the embedding method. We plan to discuss this in more detail in \cite{pwy}.

{\bf Acknowledgement}. The first and second authors would like to thank Jeff Cheeger for interesting discussions and constant encouragement.

\sect{Moser's weak maximum principle}

Let $N$ be a compact $n$-dimensional manifold with
non-empty boundary, $n\ge 3$. For a given metric $g$ on $N$
the Sobolev constant $C_S$ of $g$ (or better, of the Riemannian manifold
$(N,g)$) is defined to be the supremum of
$(\int |f|^{2n\over n-2})^{n-2\over 2n}$ over all 
 $C^1$ functions $f$ on $N$ with
$\int|\nabla f|^2=1$, which vanish along $\partial N$. So for such $f$ we have
\be
\|f\|_{\frac{2n}{n-2}} \leq C_S \| \nabla f \|_2.
\ee

\begin{theo}  \label{max}
Let $f,b$ be smooth nonnegative functions on $N \times [0,T]$ which satisfy the following:
\be \label{heat}
\frac{\partial f}{\partial t} \leq \Delta f + bf \ \mbox{on} \ N \times [0,T],
\ee
where $\Delta$ is the Laplace-Beltrami operator of the metric $g(t)$, $b$ is assumed to satisfy:
\[
\sup_{0\leq t \leq T} \left(\int_N b^{q/2} \right)^{2/q} \leq \beta~,
\]
for some $q > n$. Put $l=\max\limits_{0\le t\le T}
\left\|{\partial g\over \partial t}\right\|_{C^0}$ (norm measured in
$g(t)$) and $C_S=\max\limits_{0\le t\le T} C_S(g(t))$.
 Then, given $p_0 > 1$, there exsits a constant $C = C(n,q,p_0, \beta, C_S, l, T,R)$ such that for any $x$ in the interior of $N$ and $t \in (0,T]$,
\be  \label{infe}
|f(x,t)| \leq C t^{-\frac{n+2}{2p_0}} \left( \int_0^T \int_{B_R} f^{p_0} \right) ^{1/p_0},
\ee
where $R={1\over 2} dist_{g(0)}(x,\partial N)$ and $B_R=B_R(x)$ denotes
the geodesic ball of radius $R$ and center $x$ defined in terms of the metric 
$g(0)$.
\end{theo}
\Pf
We follow closely [Ye], [Yn1].  Let $\eta$ be a non-negative Lipschitz function vanishing
along $\partial N$. The partial differential inequality
(\ref{heat}) implies for $p\geq 2$
\[
\frac{1}{p}\frac{\partial}{\partial t}\int f^p\eta^2
dv_t\leq \int\eta^2f^{p-1}\Delta f\,dv_t+\int
bf^p\eta^2dv_t+\frac{1}{p}\int f^p\eta^2\frac{\partial}{ \partial t}dv_t~,
\]
where $dv_t$ is the volume form of $g(t)$.
Integration by parts yields
\ban
\int\eta^2f^{p-1}\Delta f\,dv_t & = &
-\frac{4(p-1)}{p^2}\int |\nabla (\eta f^{p/2})|^2
dv_t+\frac{4}{ p^2}\int |\nabla \eta|^2 f^p dv_t  \\
& & +\frac{4(p-2)}{ p^2}\int \nabla (\eta f^{p/2})
f^{p/2}\nabla \eta\,dv_t \\
& \leq &   -\frac{2}{ p}\int |\nabla (\eta f^{p/2})|^2
+\frac{2}{ p}\int |\nabla \eta|^2 f^p,
\ean
where the gradient $\nabla$ refers to $g(t)$.  By H\"{o}lder inequality
\ban
 \int bf^p\eta^2dv_t & \leq & (\int b^{q/2}dv_t)^{2/q} \left(\int (f^{p}\eta^2)^{\frac{q}{q-2}}dv_t\right)^{\frac{q-2}{q}} \\
& \leq & \beta ( \epsilon^{-\frac{n-2}{q}} 
\int f^p \eta^2 dv_t )^{1-\frac{n}{q}} \left(\epsilon^{1-\frac{n}{q}}
 \int (f^{p}\eta^2)^{\frac{n}{n-2}}dv_t\right)^{\frac{n-2}{n} \frac{n}{q}} \\
& \leq & \beta \left(\epsilon^{-\frac{n-2}{q}}\int f^p \eta^2 dv_t\right)^{1-n/q} \left(\epsilon^{(1-\frac{n}{q})(\frac{n-2}{n})}C_S^2 \int |\nabla (\eta f^{p/2})|^2 dv_t\right)^{n/q} \\
& \leq & (1-\frac{n}{q}) \beta \epsilon^{-\frac{n-2}{q}}\int f^p\eta^2 dv_t+  \frac{n}{q}\beta\epsilon^{(1-\frac{n}{q})(\frac{n-2}{n})}  C_S^2 \int |\nabla (\eta f^{p/2})|^2 dv_t.
\ean
Therefore
\ban  
\lefteqn{\frac{1}{ p}\frac{\partial}{\partial t} \int f^p\eta^2 \leq 
 \left(-\frac{2}{ p} + \frac{n}{q}\epsilon^{(1-\frac{n}{q})(\frac{n-2}{n})}\beta C_S^2\right)\int |\nabla (\eta f^{p/2})|^2 } \\
& & +\frac{2}{ p}\int |\nabla \eta|^2 f^p 
+ \left(\beta (1-\frac{n}{q}) \epsilon^{-\frac{n-2}{q}} + \frac{\sqrt nl}{p}\right)\int f^p \eta^2~.
\ean
(The volume form is omitted.)

Setting 
\[
\epsilon = (\frac{q}{np\beta C_S^2})^{\frac{nq}{(q-n)(n-2)}},
\]
we obtain the following basic estimate:
\be \label{e1}
\frac{\partial}{\partial t} \int f^p\eta^2 + \int |\nabla (\eta f^{p/2})|^2 \leq 2\int |\nabla \eta|^2 f^p + C_1(n,p, \beta, C_S,q,l) \int f^p \eta^2.
\ee

Now given $0<\tau<\tau'<T$, let
\[
\psi(t)=\left\{ \begin{array}{ll} 0&  0\leq t\leq\tau~,\\
(t-\tau)/(\tau'-\tau)&  \tau\leq t\leq \tau'~,\\
1&  \tau'\leq t\leq T~. \end{array} \right.
\]
Multiplying (\ref{e1}) by $\psi$, we obtain
\[
\frac{\partial}{\partial t}\left(\psi\int f^p\eta^2\right)+\psi
\int|\nabla (\eta f^{p/2})|^2\leq 2\psi\int|\nabla\eta|^2 f^p
+(C_1\psi+\psi') \int f^p\eta^2~.
\]
Integrating this with respect to $t$ we get
\[
\int_t f^p\eta^2+\int^t_{\tau'}\int|\nabla (\eta f^{p/2})|^2
\leq 2\int^T_\tau \int |\nabla\eta|^2 f^p
+\left( C_1+\frac{1}{ \tau'-\tau}\right)
\int^T_\tau \int f^p \eta^2~, \ \ \tau' \leq t \leq T.
\]
Applying this estimate and the Sobolev inequality we deduce
\ba  \label{e2}
\int^T_{\tau'}\int f^{p(1+\frac{2}{ n})}\eta^{2+\frac{1}{ n}}
& \leq & \int^T_{\tau'}
\left(\int f^p\eta^2\right)^{2/n}\left(\int f^{\frac{pn}{ n-2}}
\eta^{\frac{2n}{n-2}}\right)^{\frac{n-2}{n}} \nonumber\\
&\leq &  C_S^2\Bigl(\sup_{\tau'\le t\le T} \int f^p\eta^2\Bigr)^{2/n}
\int^T_{\tau'} \int \left|\nabla(\eta f^{p/2})\right|^2  \nonumber\\
&\leq &  C_S^2\left[ 2\int^T_\tau \int
|\nabla\eta|^2 f^p+\left(C_1+\frac{1}{\tau'-\tau}\right)
\int^T_\tau \int f^p \eta^2\right]^{1+\frac{2}{n}}~.
\ea

For a given $x$ in the interior of $N$, $p \geq p_0$ and $0 \leq \tau \leq T$, we put
\[
H(p,\tau,R)=\int^T_\tau\int_{B_R} f^p~,
\]
where $B_R$ is the geodesic ball of center $x$ and radius
$R$ measured in $g(0)$.
Choosing a suitable cut-off function $\eta$ ( and noticing 
$|\nabla \eta|_t\le |\nabla \eta|_0 e^{{1\over 2}lt}$)
we derive from (\ref{e2})
\be  \label{e4}
H\left(p\left(1+{2\over n}\right),\tau',R\right) \leq C_S^2
\left[ C_1+\frac{1}{ \tau'-\tau}+
{2 e^{lT}\over (R'-R)^2}\right]^{1+\frac 2n}
H(p,\tau,R')^{1+\frac 2n}~,
\ee
where $0<R<R'\le dist_{g(0)}(x,\partial N)$. To proceed we set
$\mu=1+{2\over n},\; p_k=p_0\mu^k,\; \tau_k=(1- {1\over \mu^{k+1}})t$
and $R_k={R\over 2}(1+{1\over \mu^{k/2}})$ with $R={1\over 2} dist_{g(0)}
(x,\partial N)$. Then it follows from (\ref{e4}) that
\ban
\lefteqn{H(p_{k+1},\tau_{k+1}, R_{k+1})^{1/p_{k+1}}\leq} \\
& &  C_S^{2/p_{k+1}}\left[C_1+{\mu^2\over \mu-1}\cdot {1\over 
t}+
\frac{4 e^{lT}}{R^2}\cdot {\mu\over (\sqrt \mu-1)^2}\right]^{1/p_{k}}
\mu^{k/p_k}H(p_k,\tau_k,R_k)^{1/p_k}~.
\ean
Hence
\ban
\lefteqn{ H(p_{m+1},\tau_{m+1},R_{m+1})^{1/p_{m+1}}
\leq} \\
& &  C_S^{\sum^m_{k=0}{2\over p_{k+1}}}
\left[C_1+{\mu^2\over \mu-1}\cdot{1\over t}+
\quad {4 e^{lT}\over R^2}\cdot {\mu\over (\sqrt\mu-1)^2}\right]^
{\sum^m_{k=0}{1\over p_k}}
\mu^{\sum^m_{k=0}{k\over p_k}}
 H(p_0,\tau_0, R_0)^{1/p_0}~.
\ean 
Letting $m\to\infty$
we conclude
$$|f(x,t)|\leq (1+\frac{2}{n})^{\frac{n^2+2n}{4p_0}} C_S^{\frac{n}{p_0}}\left(C_1+{1\over t}+
{e^{lT}\over R^2}\right) ^{n+2\over 2p_0}
\Bigl(\int^T_0\int\limits_{B_R} f^{p_0}\Bigr)^{1/p_0}~.$$
\qed

\sect{Bounding curvature tensor}

Let $(M,g_0) \in {\cal M}(n, r_0)$ for some $n, r_0$. Consider the Ricci flow
on $M$ with initial metric $g_0$. It is well-known that  the Riemann curvature
tensor $Rm$ satisfies the following evolution equation along the Ricci
flow~\cite{h}: \be \label{curv} \frac{\partial Rm}{\partial t} = \Delta Rm +
Q(Rm), \ee
where $Q(Rm)$ is a tensor that is quadratic in $Rm$. From this it follows that
\be \label{curv1}
 \frac{\partial}{\partial t}  |Rm| \leq \Delta |Rm| + c(n) |Rm|^{2}.
\ee
In order to apply Moser's weak maximum principle we need bounds on 
\[ \begin{array}{l}
\max_{0\leq t\leq T} \|Rm\|_{p_0},  \\
\max_{0\leq t\leq T} C_S, \\
\max_{0\leq t\leq T} |Ric|. 
\end{array} \]
We are going to employ suitable evolution equations associated with the Ricci
flow such as (3.1) and (3.2) to derive these bounds from the
initial bounds (i.e. bounds at $t=0$). The basic logic of the argument goes as
follows. By smoothness of the Ricci flow, the initial bounds can at least be
extended to a very short time interval. We then derive a uniform estimate for
the length of the maximal interval of extension.   

To get initial bounds we
actually have to pull the metric up to the tangent spaces. For any fixed $x \in
M$, we lift the metric $g_0$ by the exponential map to $\tilde{g}_x(0)$ on
$\tilde{B}_{r_0} (x) \subset T_xM$. Then \[ | \Ric_{\tilde{B_{r_0}}}| \leq 1, \ \
inj_{\tilde{B}_{r_0/2}} \geq r_0/2. \]
By \cite{a2}, the metric $\tilde{g}_x(0)$ also comes with the following uniform bounds:
\ban
\|Rm\|_{p_0, \tilde{B}_{r_0/2}(x)} & \leq & K(n,r_0,p_0), \ \mbox{for all}\ 1 \leq p_0 < \infty, \\
C_S (\tilde{B}_{r_0/2}(x)) &\leq &\chi (n,r_0).
\ean
This furnishes us with the required initial bounds.

Let $\tilde{g}_x(t)$ be the lifted Ricci flow on $\tilde{B}_{r_0/2} (x)$.
Choose $p_0 > n/2+1$, e.g. $p_0 = n+2$. Then we have

  \begin{prop} \label{contbd}
Let $[0,T_{\max} )$ be a maximal time interval on which the Ricci flow $g(t)$
 has a smooth solution such that each lifted Ricci flow $\tilde{g}_x(t)$ on
$\tilde{B}_{r_0/2}(x)$ satisfys the following estimates: \ba
\|Rm (\tilde{g}_x (t))\|_{p_0,\tilde{B}_{r_0/4}(x)}& \leq &2K(n,r_0,p_0), 
\ \ \  \label{l2} \\
|\Ric (\tilde{g}_x(t))| & \leq&  2, \label{l3} \\
C_S (\tilde{g}_x (t)) & \leq & 2\chi(n,r_0). \label{l1} 
\ea
Then $T_{\max} \geq T(n, r_0,p_0)$ for some $T(n, r_0, p_0) > 0$. 
\end{prop}

The proof of this proposition requires four lemmas. Set $K=K(n, r_0, p_0)$.

\begin{lem}
For $0 \leq t < T_{\max}$ there holds
\be  \label{r1}
\|Rm (\tilde{g}_x (t))\|_\infty \leq C(n,p_0,\chi,K) t^{-\frac{n+2}{2p_0}}.
\ee
\end{lem}

This follows immediately from Moser's weak maximum principle Theorem~\ref{max}. 

\begin{lem} \label{lem1}
For $0 \leq t < T_{\max}$ there holds 
\be  \label{me}
\max_{x \in M} \|Rm (\tilde{g}_x (t))\|_{p_0,\tilde{B}_{r_0/4}(x)} \leq \frac{K}{1-C(n,p_0,\chi,K)t}.
\ee
\end{lem}

\begin{lem} \label{lem2}
For $0 \leq t < T_{\max}$ there holds 
\be  \label{m1}
|\Ric (g(t))| \leq  e^{C(n,p_0,K,\chi)t^{\frac{2p_0-n-2}{2p_0}}}.
\ee
\end{lem}

\begin{lem} \label{lem3}
For $0 \leq t < T_{\max}$ there holds
\be  \label{m2}
\max_{x \in M}C_S (\tilde{g}_x(t)) \leq \chi e^{C(n,p_0,K,\chi)t^{\frac{2p_0-n-2}{2p_0}}}.
\ee
\end{lem}

The proof of Lemmas~\ref{lem1}-\ref{lem3} is presented in the next section.
Granting the lemmas the proof of Proposition~\ref{contbd} can be seen as
 follows. First of all we quote the following theorem of Hamilton \cite{h}:
\begin{theo}
Let $g(t), \ 0 \leq t < T$ be a (smooth) solution to the Ricci flow (\ref{rf}) on a compact manifold $M$. If the sectional curvatures of $g(t)$ remain bounded as $t \ra T$, then the solution extends smoothly beyond $T$.
\end{theo}
By (\ref{me}), (\ref{m1}), (\ref{m2}) there exists $T(n, r_0,p_0) >0$ such 
that if $t< \min(T(n, r_0,p_0), T_{\max})$ then strict inequalities hold
in~(\ref{l2}), (\ref{l3}), (\ref{l1}). Consequently, if $T_{\max} < T(n,
r_0,p_0)$, by (\ref{r1}) and Hamilton's theorem, the solution to the Ricci flow
can be extended beyond $T_{\max}$ with~(\ref{l2}), (\ref{l3}), (\ref{l1}) still
holding, contradicting the maximality of $T_{\max}$. \newline

\noindent {\em Proof of Theorem~\ref{main}}. Theorem~\ref{main} follows immediately from Proposition~\ref{contbd} and 
(\ref{r1}) if we let $p_0 = n+2$.
\qed

\sect{Proof of the three basic lemmas}

{\em Proof of Lemma~\ref{lem1}}. By (\ref{curv}) 
 and the bounds implied by the definition of $T_{max}$, we can apply (\ref{e1}) with $p=p_0$ to deduce
  \be
\frac{\partial}{\partial t} \int_{\tilde{B}_{r_0/2}(x)} |Rm|^p\eta^2 \leq 2\int_{\tilde{B}_{r_0/2}(x)} |\nabla \eta|^2 |Rm|^p + C_1(n,p_0, K, \chi) \int_{\tilde{B}_{r_0/2}(x)} |Rm|^p \eta^2.
\ee
We choose $\eta$ so that the following is true.
\be  \label{ee}
\frac{\partial}{\partial t} \int_{\tilde{B}_{r_0/4}(x)} |Rm|^p \leq (2+C_1)\int_{\tilde{B}_{r_0/2}(x)}|Rm|^p.
\ee
Now cover the ball $\tilde{B}_{r_0/2}(x)$ by balls $\tilde{B}_{r_0/4}(\tilde{y}_i)$ with $\tilde{y}_i \in \tilde{B}_{r_0/2}(x)$. By the following proposition of Gromov the number of covering is uniformly bounded.

\begin{prop}[Gromov] (see \cite[Proposition 3.11]{c}) \label{bc}
Let the Ricci curvature of $M^n$ satisfy $\Ric_{M^n} \geq (n-1)H$. Then given $r,\epsilon > 0$ and $p \in M^n$, there exists a covering, $B_p(r) \subset \cup_1^N B_{p_i} (\epsilon), \ (p_i$ in $B_p(r))$ with $N \leq N_1 (n, Hr^2, r/\epsilon)$. Moreover, the multiplicity of this covering is at most $N_2(n,Hr^2)$.
\end{prop}

Therefore $\tilde{B}_{r_0/2}(x) \subset \bigcup_1^N \tilde{B}_{r_0/4}(\tilde{y}_i)$ with $N \leq N(n,r_0)$. 

Let $y_i = \exp_x \tilde{y_i}$. We now construct a map
\be
\exp^{-1}_{y_i} \circ \exp_x: \ \tilde{B}_{r_0/4}(\tilde{y}_i) \subset T_xM \ra \tilde{B}_{r_0/4}(y_i) \subset T_{y_i}M.
\ee
For any $\tilde{y}\in \tilde{B}_{r_0/4}(\tilde{y}_i)$ connect it with the center $\tilde{y}_i$ by the unique minimal geodesic. This projects down to a geodesic on $M$ which can then be lifted to a geodesic in $\tilde{B}_{r_0/4}(y_i)$. The end point of this geodesic gives the image of $\tilde{y}$.  Thus defined,  $ \exp^{-1}_{y_i}\circ \exp_x$ is one to one and an isometry. It follows then
\[
\int_{ \tilde{B}_{r_0/4}(\tilde{y}_i)} |Rm|^p= \int_{\tilde{B}_{r_0/4}(y_i)}|Rm|^p.
\]
This combines with (\ref{ee}) implies
\[
\frac{\partial}{\partial t} \int_{\tilde{B}_{r_0/4}(x)} |Rm|^p \leq 
(2+C_1) N \max_{y \in M} \int_{\tilde{B}_{r_0/4}(y)} |Rm|^p.
\]
Integrating this, we obtain
\[
\int_t |Rm(\tilde{g}_x)|^p \leq \int_{t=0} |Rm(\tilde{g}_x)|^p + (2+C_1) N t \max_{y \in M}\int |Rm(\tilde{g}_y)|^p.\]
Therefore
\[
(1-(2+C_1) N t)\max_{y \in M}\int |Rm(\tilde{g}_y)|^p \leq \int_{t=0} |Rm(\tilde{g}_x)|^p \leq K,
\]
and this implies (\ref{me}).
\qed
\newline

{\em Proof of Lemma~\ref{lem2}}. Note that for Ricci curvature we do the estimate directly, without passing to the tangent space.

 Recall that the Ricci curvature tensor satisfies the following evolution equation \cite{h}.
\[
\frac{\partial}{\partial t} \Ric = \Delta \Ric + 2Rm ( \Ric) - 2\Ric \cdot \Ric.
\]
Let $\varphi (t) = \max_{x \in M} |\Ric (g(t))|$, then
\[
\frac{\partial}{\partial t} \varphi (t) \leq c(n) |Rm (g(t))| \varphi (t).
\]
Using (\ref{r1}), together with $\varphi (0) \leq 1$ and integrating with respect to $t$ gives
\[
\varphi (t) \leq e^{\int_0^t C(n,p_0,K,\chi) \tau^{-\frac{n+2}{2p_0}}d\tau} = 
e^{ \frac{2p_0}{2p_0-n-2} C(n,p_0,K,\chi) t^{\frac{2p_0-n-2}{2p_0}}}.
\]
This implies (\ref{m1}).
\qed
\newline

{\em Proof of Lemma~\ref{lem3}}.
For $u \in C^\infty(\tilde{B}_{r_0/2} (x))$ and vanishing on the boundary, define
\[
E_t[u] =\left( \frac{\|u\|_{\frac{2n}{n-2}, \tilde{B}_{r_0/2}}}{\|\nabla u\|_2}\right)^2.
\]
Then a straightforward computation shows that
\[
\frac{\partial}{\partial t} E_t [u] \leq c(n) \|\Ric (\tilde{g}_x(t))\|_\infty E_t[u].
\]
Integrating this gives
\[
C_S (\tilde{g}_x(t)) \leq \chi e^{C(n,p_0,K,\chi)t^{\frac{2p_0-n-2}{2p_0}}},
\]
which is (\ref{m2}).
\qed

\sect{Applications}
By Theorem~\ref{main}, the deformed metric $g(t)$ has uniform sectional
 curvature bound (away from $t=0$) and $g(t)$ is close to $g(0)$. To apply the
results with sectional curvature bounds to $g(t)$, we need to show that other
geometric quantities like diameter, volume and excess are also under control. We
first prove the following lemma. 
\begin{lem} Let $g(t)$ be the Ricci flow in
Theorem~\ref{main}. Then for $0 \leq t \leq T(n,r_0)$,
 \ba e^{-2t} \diam (0)
\leq & \diam (t) & \leq  e^{2t} \diam (0), \label{dm}\\
 e^{-4nt} \vol (0) \leq &
\vol (t) & \leq e^{4nt} \vol (0), \label{vm}\\ & ex (t) & \leq e^{2t} ex
(0)+(e^{2t}-e^{-2t})\diam (0), \label{exs} \ea
where e.g. \diam (t) means \diam(g(t)).
\end{lem}
\Pf For the estimate on the diameter we consider the evolution of the length functional. Thus fix a curve $c$ and let $l_c(t)$ denote its length in the metric $g(t)$. Then
\be 
-2 l_c(t) \leq \frac{\partial l_c(t)}{\partial t} \leq 2 l_c(t), \label{vl}
\ee
from which we obtain $e^{-2t}l_c(0) \leq l_c(t) \leq e^{2t}l_c(0)$. This gives
\be \label{dl}
e^{-2t} d_{p,g} (0) \leq d_{p,q}(t) \leq e^{2t} d_{p,g} (0),
\ee
where $d_{p,q}(t)$ denote the distance between $p$ and $q$ in the metric $g(t)$ for $p,\ q$ in $M$, in particularly (\ref{dm}).

For the volume let $\omega(t)$ denote the volume form (or density if $M$ is not orientable) of $g(t)$. Consider $A(t)=\omega(t)/\omega(0)$. One computes
\[  \frac{\partial A(t)}{\partial t}=\Tr_{g(t)} \frac{\partial g(t)}{\partial t}) A(t).
\]
Then
 \[
-4n A(t) \leq \frac{\partial A(t)}{\partial t} \leq 4n A(t).\]
Hence $e^{-4nt} \leq A(t)\leq e^{4nt}$ and the volume estimate follows.

(\ref{exs}) follows from (\ref{dl}).          \qed 
\newline

\noindent {\em Proof of Theorem~\ref{aflat}}. Let $g(t)$ be the unique solution to (\ref{rf}) with the given metric as the initial data. By (\ref{t2}), for $0 < t \leq T(n,r_0)$,
\be  \label{kf}
|K(t) |  \leq  C(n,r_0) t^{-1/2}
\ee
If $\epsilon_0(n)$ is the small constant in the original Gromov's almost flat manifold theorem \cite{g1}, we choose 
 $t_0 = T(n,r_0), \ \epsilon = \diam (0) \leq\left(\frac{\epsilon_0 (n)t_0^{1/2}}{C(n,r_0)e^{4t_0}}\right)^{1/2}$. Then from (\ref{dm}) and (\ref{kf}) $|K(t_0) D^2(t_0)| \leq \epsilon_0 (n)$. Applying Gromov's almost flat manifold theorem to $g(t_0)$ gives Theorem~\ref{aflat}.
\qed
\newline

\noindent {\em Proof of Theorem~\ref{fib}}. Consider the Ricci flows $g_M(t), \ g_N(t)$  on $M$ and $N$ respectively, starting with the given metrics. Then
\[
|K_M(t)| \leq C(n,r_0) t^{-1/2}, \ \ |K_N(t)| \leq C(m,\mu) t^{-1/2}, 
\]
for $0 < t \leq T(n,r_0)$ and $0 < t \leq T(m,\mu)$ respectively. Let 
$T=\min(T(n,r_0),\, T(m,\mu))$, and $C=\max( C(n,r_0),\, C(m,\mu))$. Rescale the metrics
\[
h_M(t) =  C t^{-1/2} g_M(t), \ \ h_N(t) =  C t^{-1/2} g_N(t)
\]
so that with respect to $h_M(t), h_N(t)$, where $0<t\leq T$,
\be
|K_M| \leq 1,  \ \ |K_N| \leq 1.
\ee

Since $\inj_N \geq \mu$ for the initial metric on $N$, we have $C_S (g_N(0)) \leq \chi (m, \mu)$ and therefore by (\ref{l1})
\[
C_S (h_N(t)) = C_S(g_N(t)) \leq 2 C_S (g_N(0)) \leq 2\chi (m, \mu),
\]
for $0 < t \leq T$. By \cite{cgt} 
\[
\inj_{h_N(t)} \geq \mu_1 (m, \mu)
\]
for all $0 < t \leq T$.

Now for $0 < t \leq T$,
\ban
\lefteqn{d_H\left(M(h(t)), N(h(t))\right) = C t^{-1/2} d_H (M(g(t)), N(g(t))) }\\
&\leq &  C t^{-1/2} \left[ d_H (M(g(t)), M(g(0))) + d_H (M(g(0)), N(g(0))) + d_H (N(g(0)), N(g(t))) \right] \\
& \leq & C t^{-1/2} \left[ 8t + d_H (M(g(0)), N(g(0))) \right].
\ean
Let $\lambda (n)$ be the small constant in the refined fibration theorem (cf. \cite[Theorem 2.6]{cfg}). We choose $t_0 \leq \min \{ T, \left( \frac{\mu_1(n,\mu) \lambda (n)}{16 C} \right)^2\}, \  \epsilon =d_H (M(g(0)), N(g(0))) \leq \frac{1}{2} \left( \frac{\mu_1(n,\mu) \lambda (n) t_0^{1/2}}{C} \right)$. Then $\frac{d_H\left(M(h(t_0)), N(h(t_0))\right)}{\inj_{h_N(t_0)}} \leq \lambda (n)$. Applying Theorem 2.6 in \cite{cfg} finishes the proof.
\qed
\newline

With (\ref{dm}) and (\ref{exs}) the proof of Theorem~\ref{soul} is quite similar. We observe that the original statement in \cite[Theorem 2.1]{pz} can be restated as that there exists an $\epsilon (n)$ such that if $M^n$ is a compact manifold satisfying $| K_M| \leq 1, \ ex(M)/diam(M) \leq \epsilon (n)$, then either $M$ is infranil or diffeomorphic to the union of two normal bundles over two embedded infranilmanfolds in $M$. We then follow the above argument.

It should be pointed out that Theorem~\ref{soul} implies a sphere theorem by some purely topological argument, see \cite{pz} for details.

Theorem~\ref{minv} follows from Theorem 1.1, (\ref{vm}) and Corollary 0.4 in
\cite{ro1}. 

Theorem~\ref{h2} and \ref{fg} can be proved by using Theorem~\ref{fib} as in \cite{ro2}.


X. Dai

Department of Mathematics, 
University of Southern  California, 
Los Angles, CA  90089 

 xdai@math.usc.edu

G. Wei \& R. Ye

Department of Mathematics,
 University of California,
Santa Barbara, CA 93106

 wei@math.ucsb.edu \ \ \ 
yer@math.ucsb.edu

\end{document}